\documentclass[10pts]{article}
\usepackage{latexsym}
\usepackage{amsmath}
\usepackage{amsfonts}
\usepackage{epsfig}
\usepackage{amssymb}
\usepackage{amsthm}
\usepackage{shuffle}

\input xy
\xyoption{all}

\title{A few conjectures about the multiple zeta values}
\author{German Combariza}

\newtheorem{Def}{Definition}
\newtheorem{Teo}[Def]{Theorem}
\newtheorem{Lema}[Def]{Lemma}

\newtheorem{Cor}[Def]{Corollary}
\newtheorem{Con}[Def]{Conjecture}

\theoremstyle{definition}
\newtheorem{Exa}[Def]{Example}

\newcommand{\R}{\mathbb{R}}

\newcommand{\Q}{\mathbb{Q}}

\newcommand{\ses}{s_1,\cdots,s_k}
\newcommand{\sk}{s_1,s_2,\cdots,s_k}
\newcommand{\xk}{x^{s_1-1}yx^{s_2-1}y\cdots x^{s_k-1}y}

\begin{document}
\maketitle

\begin{abstract}
The multiple zeta values (MZV) are a set of real numbers with a beautiful structure as an algebra over the rational numbers. They are related to maybe the most important conjecture on mathematics today, the Riemann hypothesis. In this paper we will show partial solutions to five well known conjectures about the MZV which we partially proved using just linear systems and shuffles products. We saw that to partially solve four of this conjectures we only need linear algebra.
\end{abstract}

\section{Preliminary}
Let $s=(s_1,s_2,\cdots, s_k)$ be a tuple of natural numbers with $s_i>0$ and $s_1>1$. Consider the real numbers also known as the \textit{multiple zeta values} (MZV)
$$
\zeta(s) = \sum_{n_1>n_2>\cdots>n_k>0} \frac{1}{n_1^{s_1}n_2^{s_2}\cdots n_k^{s_k}}.
$$
Our main interest is the problem concerning the polynomial relations over $\Q$ of the MZV. For MZV of length 1 and even degree there is a beautiful formula due to Euler
$$
\zeta(s)=-\frac{(2\pi i)^sB_s}{2s!}\,\,\,\text{ for } s=2,4,6,\cdots
$$
which express the values of the zeta functions at even points in terms of the Bernoulli numbers $B_s\in\Q$ defined by the generating function
$$
\frac{t}{e^t-1} = 1 -\frac{t}{2} + \sum_{s=2}^\infty B_s \frac{t^s}{s!}\,\,\,\, B_s=0\,\,\,\,\text{ for odd } s\ge3.
$$
Therefore the transcendence degree of the ring $\Q[\zeta(2),\zeta(4),\cdots]$ over $\Q$ is 1. Much less is known on the arithmetic nature of values of the zeta function at odd integers $\zeta(3), \zeta(5), \cdots$. We believe that each one of these numbers is transcendent, we only know by Apery's theorem that $\zeta(3)$ is an irrational number. On the other hand for the multi-index MZV there are more relations explained by a rich algebraic structure over an associate algebra which surjects the algebra of the MZV over $\Q$. For example for the MZV $\zeta(2,1)$ Euler  proved the identity 
$$
\zeta(2,1)=\zeta(3).
$$
The big problem is: \textit{Find all the relations between the MZV}. Apparently all the relations are explained by the shuffle and stuffle products and maybe the main conjecture about it is the Zagier's conjecture which determine the dimension of the $\Q$-algebra spam by the MZV of a fixed degree over $\Q$. We will start this paper by defining the shuffle and stuffle products in section 2. We will also state the Zagier's conjecture. The goal of section 3 is explain the Kaneko-Noro-Tsurumaki conjecture \cite{Kan Nor Tsu}. In section 4 we will explain the Petitot-Minh approach to the Zagier's conjecture using the xtaylor algorithm. In section 5 we will talk about the transcendence degree over $\Q$ and the 2-3 conjecture. Finally in section 6 the Broadhurst-Kreimer conjecture and decomposability of the MZV. We started this project inspired by the Petitot-Minh paper \cite{Min Pet} and we will follow its notation and definitions.

\section{A source of relations}

Consider the alphabet $X={x_0,x_1}$ and encode the multi-index $s=(\ses)$ by the rule
\begin{eqnarray}\label{encode}
s=(\ses)\mapsto x_0^{s_1-1}x_1x_0^{s_2-1}x_1\cdots x_0^{s_k-1}x_1.
\end{eqnarray}

Then the degree of the MZV $\zeta(s)$ is by definition the degree of the monomial on the variables $x_0,x_1$, i.e. $s_1+s_2+\cdots s_k$. Let $\mathcal{H}^1$ and $\mathcal{H}^2$, be the sub-algebras $\Q\oplus\Q\langle X\rangle x_1$ and $\Q\oplus x_0\Q\langle X\rangle x_1$, respectively.  Let us define the the shuffle product $\shuffle$ and the stuffle or second shuffle product $*$ over $Q\langle X\rangle$ by the rules
\begin{eqnarray*}
    1\shuffle w = w\shuffle 1 = w && 1*w=w*1=w.
\end{eqnarray*}
and 
\begin{eqnarray}\label{shuffle and stuffle def}
    x_iu\shuffle x_jv& = &x_i(u\shuffle x_jv) + x_j(x_iu\shuffle v)\\
    y_iu*y_jv& = &y_i(u* y_jv) + y_j(y_iu* v) + y_{i+j}(u* v)
\end{eqnarray}
Inductive arguments enable us to prove that each of the above products is commutative and associative. The motivation for this products are explained in the following examples.
\begin{Exa}[The Stuffle product]
Using $\Sigma$ notation, the MZV $\zeta(2)^2$ can be expressed as
$$
\zeta(2)^2 = \left(\sum_{m>0}\frac{1}{m^2}\right)\left(\sum_{n>0}\frac{1}{n^2}\right)=\left(\sum_{m>n>0}+\sum_{n>m>0}+\sum_{m=n}\right)\frac{1}{m^2n^2} 
$$
This is
$$
   \zeta(2)^2 = 2\zeta(2,2) + \zeta(4).
$$
Using stuffle notation, should be clear that
        $$y_2*y_2 = y_2(1*y_2) + y_2(y_2*1) + y_4 = 2y_2y_2 + y_4.$$
\end{Exa}

The MZV can also be expressed as multiple integrals.
    $$ \zeta(\sk) = \int\cdots\int_{1>t_1>\cdots>t_p>0} w_1(t_1)w_2(t_2)\cdots w_p(t_p). $$
    with $p=s_1+s_2+\cdots s_k$ and $w_i(t) = dt/(1-t)$ if $i\in\{s_1,s_1+s_2,\cdots,s_1+s_2+\cdots+s_k\}$ and $w_i(t)=dt/t$ otherwise. 
\begin{Exa}[The Shuffle product] Using the integral notation above for the MZV $\zeta(2)^2$ we have
$$\zeta(2)^2 = \left(\int_{0>t_1>t_2>1}\frac{dt_1dt_2}{t_1(1-t_2)} \right)\left(\int_{0>r_1>r_2>1}\frac{dr_1dr_2}{r_1(1-r_2)} \right)$$
simplifying
$$4\left(\int_{0>t_1>t_2>t_3>t_4>1}\frac{dt_1dt_2dt_3dt_4}{t_1t_2(1-t_3)(1-t_4)}\right)+
  2\left(\int_{0>t_1>t_2>t_3>t_4>1}\frac{dt_1dt_2dt_3dt_4}{t_1t_3(1-t_2)(1-t_4)} \right)$$
i.e.
$$\zeta(2)^2 = 4\zeta(3,1) + 2\zeta(2,2).$$
Using the shuffle product the same identity is
     \begin{eqnarray*}
         x_0x_1\shuffle x_0x_1 &=& x_0(x_1\shuffle x_0x_1) + x_0(x_0x_1\shuffle x_1) = 2x_0(x_1\shuffle x_0x_1) \\
                       &=& 2x_0[x_1(1\shuffle x_0x_1) + x_0(x_1\shuffle x_1)] = 2x_0x_1x_0x_1 + 4x_0^2x_1^2 \\
            \zeta(2)^2 &=& 2\zeta(2,2) +4\zeta(3,1).
	\end{eqnarray*}
\end{Exa}

Consider the function $\zeta$ as the $\Q$-linear map  $\zeta:\mathcal{H}^2\to\R$ induce by:
$$
\zeta:\xk\mapsto\zeta(s) = \sum_{n_1>n_2>\cdots>n_k>0} \frac{1}{n_1^{s_1}n_2^{s_2}\cdots n_k^{s_k}}.
$$
\begin{Teo}
Under the shuffle product, the $\zeta$ map $\mathcal{H}^2\to \R$ is an homomorphism.
$$ \zeta(w\shuffle v)=\zeta(w)\zeta(v) $$
\end{Teo}
The same is true for the stuffle product. For the proof of this theorems we refer the reader to \cite{Min Pet}
\begin{Teo}
Under the stuffle product, the $\zeta$ map $\mathcal{H}^2\to \R$ is an homomorphism.
$$ \zeta(w* v)=\zeta(w)\zeta(v) $$
\end{Teo}
\begin{Cor} For any two words $w, w'$ in $\mathcal{H}^2$
$$
\zeta(w\shuffle w' - w*w')=0.
$$
\end{Cor}
This beautiful corollary is the main source of relations for the multiple zeta values and our object of study is the ideal generated by this difference.
Denote $\zeta_p$ the $\Q$-algebra form by all the MZV of weight $p$, i.e. $\zeta_p = \Q[\zeta(s)]$ for all $s$ with $s_1+\cdots+s_k=n$. Denote $\zeta_p^+$ the $\Q$-algebra spam by all the the MZV of weight less or equal $p$, i.e $\zeta_p^+ = \zeta_p\oplus\zeta_{p-1}\oplus\cdots\oplus\zeta_2$.
\begin{Def}
	The weight of the  MZV $\zeta(\sk)$ is the natural number $s_1+s_2+\cdots s_k$. 
\end{Def}
\begin{Con}[Zagier]
    	Let $d_p$ be the $\Q$-dimension of $\zeta_p$. Then these are given by the recurrence $d_1=0, d_2=d_3=1$ and 
    	$$ d_n = d_{n-2} + d_{n-3} \,\,\,\text{ for } n\ge4. $$
\end{Con}

Here are some examples of this conjecture.

\begin{itemize}
    \item $p = 2, d_2 = 1, \zeta_2 = \Q[\zeta(2)]$.
    \item $p = 3, d_3 = 1, \zeta_3 = \Q[\zeta(3)]$.
    \item $p = 4, d_4 = 1, \zeta_4 = \Q[\zeta(2)^2]$.
    \item $p = 5, d_5 = 2, \zeta_5 = \Q[\zeta(2)\zeta(3), \zeta(5)].$ 
    \item $p = 6, d_6 = 2, \zeta_6 = \Q[\zeta(2)^3, \zeta(3)^2].$ 
    \item $p = 7, d_7 = 3, \zeta_7 = \Q[\zeta(7), \zeta(2)\zeta(5), \zeta(2)^2\zeta(3)].$ 
    \item $p = 8, d_8 = 4, \zeta_8 = \Q[\zeta(2)^4, \zeta(3)\zeta(5), \zeta(2)\zeta(3)^2, \zeta(6,2)].$ 
    \item $p = 9, d_9 = 5, \zeta_9 = \Q[\zeta(9), \zeta(2)\zeta(7), \zeta(2)^2\zeta(5), \zeta(2)^3\zeta(3), \zeta(3)^3].$
\end{itemize}

\newpage
\section{Kaneko-Noro-Tsurumaki conjecture and the regularization process}
So far we have a nice way to generate relations between the MZV with the ``shuffle - stuttle'' process. The question now is the domain for this relations. One can try
$$
w\shuffle w' -w*w'
$$
for every $w,w' \in \mathcal{H}^2$ but this is going to produce more variables than relations. To solve this problem Kaneko-Noro-Tsurumaki introduce the notion of regularization. Below is a brief introduction to the subject. To view this in detail see \cite{Kan Nor Tsu}.

It is known that $(\mathcal{H}^1, \shuffle)$ is a commutative algebra and isomorphic to the polynomial algebra over $\mathcal{H}^2$ with one variable, i.e.
    \begin{center}
    $\mathcal{H}^1\simeq \mathcal{H}^2[x_1].$
    \end{center} 

Let $reg:\mathcal{H}^1\to\mathcal{H}^2$ be the map ``constant term'' from the above isomorphism of shuffle algebras.
    \begin{eqnarray}\label{regularization map}
        reg : \sum_{i=0}^n w_i\shuffle x_1^{\shuffle i}\mapsto w_0.
    \end{eqnarray} 
    
    \begin{Exa}\begin{eqnarray*}
   	\mathcal{H}^1  \simeq &\mathcal{H}^2[x_1]   &\mapsto  \mathcal{H}^2\\
        x_1x_0x_1  = &-2x_0x_1^2 + x_0x_1\shuffle x_1 &\mapsto -2x_0x_1^2
    \end{eqnarray*}\end{Exa}

This process will reduce the number of variables to elements only in $\mathcal{H}^2$ after the next theorem.

\begin{Teo}[From \cite{Kan Nor Tsu}]
    $$ \zeta(reg(w_1\shuffle w_0 - w_1*w_0)) = 0$$
    For any $w_0\in \mathcal{H}^2$ and $w_1\in \mathcal{H}^1$.
\end{Teo}

Given an integer $n$ let $w, w'\in \mathcal{H}^1$ such that $degree(w) + degree(w')=n$. Then the regularization process is giving a linear system of relations of the MZV with variables the monomials in $\mathcal{H}^2$ of degree $n$. The number of variables is $2^{n-1}$. This works, but we have limitations of space then there is a new problem: Too many equations. This is where the Kaneko-Noro-Tsurumaki conjectures plays its role.

\begin{Con}[Kaneko, Noro, Tsurumaki]
         In order to have all MZV relations of degree n it is enough to take $w_1 \in \mathcal{H}^2$ and $w_0\in\{x_1,x_0x_1, x_0^2x_1, x_0x_1^2\}$ with $degree(w_0) + degree(w_1)=n$.
\end{Con}
For $n>7$ the Kaneko-Noro-Tsurumaki linear system has $2^{n-1}$ relations and $2^{n-1}$ equations which null space is generated by the MZV relations of degree $n$. We verify this dimension with the Zagier conjecture. The files and the process on this conjecture is at the end of this paper.

\newpage
\section{Petitot-Minh approach and Zagier conjecture}
In this section we will try to explain how Petitot and Minh approached to the Zagier's conjecture using the xtaylor algorithms. We will also discuss their conjecture about the algebraic structure of the algebra of the MZV over the rationals. More details about this can be found at \cite{Min Pet}. This papers was the first motivation for this work and we follow its ideas and notation.

	\begin{Def}
		A word $l\in \Q\langle x_0,x_1\rangle$ on the letters $x_0,x_1$ is a Lyndon word if it is less that any of its proper right factors in the lexicographic order.
		$$ l=uv \Rightarrow l<v. $$
	\end{Def}
	
	The Lyndon words plays an important role in the associative algebra $(\Q\langle x_0,x_1\rangle,\shuffle)$ and determine a basis for the calculations of the MZV relations. A few examples of Lyndon are $x_0, x_1, x_0x_1, x_0^2x_1, x_0x_1^3, x_0x_1x_0x_1^2$ and examples of non Lyndon words are $x_0x_1x_0x_1, x_0x_1x_0^3x_1$. Before state the Minh-Petitot conjecture we will see that the set of Lyndon words form a basis for $\Q\langle x_0,x_1\rangle$ and we will express any words in terms of this basis.
	
	\begin{Teo}
		The set $\mathcal{L}$ of all Lyndon words over $x_0,x_1$ is a polynomial basis for the commutative shuffle algebra $\Q\langle x_0,x_1 \rangle$.
				$$ \Q\langle x,y\rangle_{\shuffle} \simeq \Q[\mathcal{L}]_{\shuffle}.$$
	\end{Teo}
	
	For the proof of this theorem and more detail we reefer the reader to \cite{Min Pet}. This the idea they followed. If one is able to express any word in $\mathcal{H}^1$ in term of Lyndon words, the variables in the systems of MZV will be reduce to the Lyndon words. To do that we need a couple of definitions that make sense the idea of a differential and the will use the same idea behind the Taylor algorithms.
	
	\begin{Def}[The bracket]
	For any Lyndon word $l$, the bracketed form $[l]$, is defined as follows:  $$[l]=[[u],[v]]$$ 
	with $[x_0]=x_0$ ,$[x_1]=x_1$ and $[x_0,x_1]=x_0x_1-x_1x_0$ and for $l=uv$, where $u,v$ are Lyndon words and $v$ is the 
	longest Lyndon word such that $l=uv$.
	\end{Def}
	\begin{Def}[The right residual]
	The right residual $(p\triangleright q)$ of $p$ by $q$ is defined by: 
	$(p\triangleright q|w)=(p| qw)$ for every word $w$. Where $(u|v)=\delta^u_v$.
	\end{Def}
	
	We are now ready to explain the xtaylor algorithms. This algorithms appeared first at Petitot's doctoral thesis. This are the steps that we will follow. First we need a derivation then find a set of possible divisor. At this stage we improve the xtaylor algorithms as appears in \cite{Min Pet} by looking to a small set of factors for $\zeta(\sk)$ to all sub-sequence of $\sk$ that are also Lyndon words i.e. $\{\zeta(s_1),\zeta(s_1,s_2,\cdots,\zeta(s_1,\cdots,s_{k-1})) \cap$ Lyndon. Once a divisor in found one proceed to do the same process for the residuals we will see below. This process follow the idea of the Taylor algorithms the only missing concept here was the differential. To those elements whose differentiation by them is not zero we are calling them here Lyndon divisors. 
	
	\begin{Lema}[2.1 \cite{Min Pet}] The operation $-\rhd[ - ]$ is a differential for the shuffle product.
	\end{Lema}
	
	We illustrate the algorithms with this example. $\zeta(2,4)$ is not a Lyndon word. The set of Lyndon divisors is just $\{\zeta(2)\}$. The Taylor process follow using $\zeta(2,4)-\zeta(2,4)\rhd[\zeta(2)]$ instead. The final stage is
	$$\zeta(2,4) = 4\zeta(5,1)+2\zeta(4,2) -\zeta(3)\shuffle\zeta(3)+\zeta(2)\shuffle\zeta(4).$$
		
	We are now ready to see the Minh-Petitot conjecture and their explanation.
	
	\begin{Lema}[4.1 \cite{Min Pet}]
		Let $R = k[X_1,X_2,\cdots,X_n]$ be a polynomial ring; $I$ be an ideal of $R$ and $G \subset R$ be a Gr\"obner. basis of $I$ for any
		 admissible order. If the leading terms of polynomials of $G$ are single indeterminates then the algebra $R/I$ is free and the ideal
		 $I$ is prime.
	\end{Lema}
	\begin{Con}
		The $\Q$-algebra of the MZV is a polynomial algebra.
	\end{Con}
	\begin{Exa}
		$\zeta_{\le 6} = \Q[\zeta(2), \zeta(3), \zeta(5)]$.
	\end{Exa}

\section{Transcendence degree of $\zeta_p$ over $\Q$ and the 2-3 conjecture}
The algebra $\zeta_6$ of MZV of degree 6 or lower over $\Q$ is generated for the MZV $\zeta(2), \zeta(3), \zeta(5)$. Therefore it's transcendence degree is 3. The natural question is to generalize this degree for higher algebras of MZV. Here another application of the Lyndon words will help to state a conjecture about it.

	\begin{Con}
		The set of MZV $\zeta(\sk)$ with $k\ge1$ and $s_j\in\{2,3\}$ for $1\le j\le k$, such that $\sk$ is a Lyndon word on the alphabet 	
		$\{2,3\}$, give a transcendence basis of the $\Q$-algebra of the MZV.
	\end{Con}
	
	If this conjecture is true we can tell what is going to be the transcendence degree if each step of this tower of $Q$-algebras.
	\begin{Cor}
	     The number $N(p)$ of elements of weight $p$ in a transcendence basis of the $MZV$ is given by
	     $$N(p)=\frac{1}{p}\sum_{l|p}\mu(p/l)P_l$$
	     where $P_l=P_{l-2}+P_{l-3}$ and $P_1=0, P_2=2, P_3=3 $.
	\end{Cor}
	
	This is the conjecture exemplified. 
	
	\begin{center}
    \begin{tabular}{ |c|c|c|} \hline
        Degree &Lyndon& MZV \\ \hline
        p = 2; &$x;$& $\zeta(2)$\\ \hline
        p = 3; &$y;$& $\zeta(3)$\\ \hline
        p = 4;&&\\ \hline
        p = 5; &$xy;$& $\zeta(2, 3)$\\ \hline
        p = 6;&&\\ \hline
        p = 7; &$x^2y;$& $\zeta(2, 2, 3)$\\ \hline
        p = 8; &$xy^2;$& $\zeta(2, 3, 3)$\\ \hline
        p = 9; &$x^3y;$& $\zeta(2, 2, 2, 3)$\\ \hline
    \end{tabular}
\end{center}

	In this context the algebra $\zeta_6$ is just the polynomial algebra $\Q[\zeta(2),\zeta(3),\zeta(2,3)]$.
	
\section{Deep vs weight and the Broadhurst-Kreimer conjecture}

This last conjecture is about the ``decomposability" of the MZV. For example $\zeta(2,3)$ can be write in term of the MZV $\zeta(2),\zeta(3)$ and $\zeta(5)$ this are MZV of deep 1 but the MZV of depth 2 $\zeta(6,2)$ can not be decompose in MZV of depth 1.

\begin{Def}
 The  depth of the MZV $\sk$ is $k$.
\end{Def}

The problem here is what and where. What MZV are not decomposable in term of  of lower depth and when can we find them.

	\begin{Con}
		The number $D_{n,k}$ of MZV of weight $n$ and depth $k$ that are not reducible to MZV of lesser depth is generated by
		$$
		1-\frac{x^3y}{1-x^2} + \frac{x^{12}y^2(1-y^2)}{(1-x^4)(1-x^6)}=\prod_{n\ge3}\prod_{k\ge1}(1-x^ny^k)^{D_{n,k}}.
		$$ 
	\end{Con}
	
	This is the table about the depth indecomposability. 
	\begin{center}
	    \begin{tabular}{ |c|c|c||c| } \hline
	    w/d & 1 & 2 &   Irr. MZV\\ \hline\hline
	    2   & 1 &   &  $\zeta(2)$ \\ \hline
		3   & 1 &   &  $\zeta(3)$\\ \hline
		4   &   &   &  \\ \hline
		5   & 1 &   &  $\zeta(5)$\\ \hline
		6   &   &   &  \\ \hline
		7   & 1 &   &  $\zeta(7)$\\ \hline
		8   &   & 1 &  $\zeta(6,2)$\\ \hline
		9   & 1 &   &  $\zeta(9)$\\ \hline
		\end{tabular}	
	\end{center}
\section{What have we done?}

\subsection{Kaneko-Noro-Tsurumaki Conjecture}. To solve the Kaneko-Noro-Tsurumaki conjecture stated in section 3 I used Gap.I generate the matrices for the MZV relations from degree 7 until degree 22. This matrices are nothing else but the difference of the shuffle and stuffle product 
$$ w_1\shuffle w_2 - w_1*w_2$$
for $w_1 \in \mathcal{H}^2$ and $w_2 \in \{x_1, x_0x_1, x_0^2x_1, x_0x_1^2\}$ with $degree(w_1) + degree(w_2) = n$. The rank of this matrices equals the the dimensions of the Zagier's conjecture. Showing that this is number of relations is enough and solving the conjecture until degree 22.

\begin{itemize}
\item Program: GAP. 
\item files in: .http://www.csd.uwo.ca/~gcombar/ResearchMZV/GAP/Regularization. 
\item function: RelationWithRegularizationTruncatedUpToDegree(  degree, filename).
\end{itemize}

\subsection{Zagier's Conjecture and Minh-Petitot Conjecture} Let us fix a degree $n$. In order to prove this conjecture for degree $n$  this are the steps I did
\begin{enumerate}
	\item Find all Lyndon words of degree $n$.
	\item Find all non-Lyndon words in $\mathcal{H}^2$ and factorize them as shuffle products of Lyndon words. For example
	$$\zeta(2,4) = 4\zeta(5,1)+2\zeta(4,2) -\zeta(3)\shuffle\zeta(3)+\zeta(2)\shuffle\zeta(4).$$
	\item I called $\zeta(2), \zeta(3)$ and $\zeta(4)$ the divisors of $\zeta(2,4)$. Here I realize that the sub-words intersect the divisors. This was a  improvement in my original xtaylor algorithm. 
	\item Find all the relations $w_1\shuffle w_2 - w_1* w_2$ for $w_1, w_2$ Lyndon words with $degree(w_1) + degree(w_2) =n$.
	\item Replace all non-Lyndon words in $w_1\shuffle w_2 - w_1* w_2$ for its correspondence factorization in Lyndon words, for example
	\item Solve the system will find a Gr\"obner. basis for the ideal generated by all the relations of degree $n$.
	\item This Gr\"obner. basis should have a linear basis of the same dimension as the Zagier's conjecture numbers.
	\item A big improvement here is use the Gr\"obner. basis for smaller degrees $\le n$ in the Lyndon factorization, for example
			$$\zeta(2,4) = 4\zeta(5,1)+2\zeta(4,2) -\zeta(3)\shuffle\zeta(3)+\zeta(2)\shuffle\left(\frac{2}{5}\zeta(2)\shuffle\zeta(2)\right).$$
\end{enumerate}
After this process find the Gr\"obner. basis is reduced to solve a linear system where the variables are products and powers of Lyndon words. Now it is easy to check the Minh-Petitot condition. 
\begin{itemize}
	\item Program: Magma.
	\item Function: Replacing( , , ).
	\item file: http://www.csd.uwo.ca/~gcombar/ResearchMZV/Magma/All/LinearSystem4.mag
\end{itemize}

\subsection{The 2-3 Conjecture}
This conjecture is easily seen from the Gr\"obner. basis and the linear generators in each degree for the algebra $\zeta_n$. For degrees $n=2,3$ the 2-3 Lyndon words and Linear generators for $\zeta_n$ are the same: $\zeta(2)$ and $\zeta(3)$ respectively. 

For $n=5$ there is one Lyndon word $\zeta(2,3)$. The linear generators for $\zeta_5$ are $\zeta(5)$ and $\zeta(2)\zeta(3)$. The connection is made by the corresponding relation in the Gr\"obner. basis of degree 5. i.e. $\zeta(2,3) = \frac{9}{2}\zeta(5) -2\zeta(2)\zeta(3)$. This shows that the transcendental degree of $\zeta_5^+$ is three with generators the three Lyndon words $\zeta(2), \zeta(3), \zeta(5)$.

\begin{eqnarray*}
    \zeta_5^+ &=&  \zeta_5 \oplus \zeta_4 \oplus \zeta_3 \oplus \zeta_2 = \zeta_5 \oplus \zeta_3 \oplus \zeta_2\\
              &=&  \Q[\zeta(5),\zeta(3)\zeta(2)]\oplus\Q[\zeta(3)]\oplus\Q[\zeta(2)] \\
              &=&  \Q[\zeta(5),\zeta(3),\zeta(2)] \\
              &=&  \Q[\zeta(2,3),\zeta(3),\zeta(2)].
\end{eqnarray*}

Just to be clear: The dimension of $\zeta_n^+$ is the transcendental degree over $\Q$. But the dimension of $\zeta_n$ is the linear dimension over $\Q$. For degree 6 there are nothing to say since there are not 2-3 Lyndon words of weight 6. \\

For degree 7 there is a 2-3 Lyndon word $\zeta(2,2,3) = - \frac{291}{16}\zeta(7) + 12\zeta(5)\zeta(2) - \frac{3}{5}\zeta(3)\zeta(2)^2$. The $\zeta_7$ space is generated by $\zeta(7), \zeta(3)\zeta(2)^2, \zeta(5)\zeta(2)$. \\

The generators for $\zeta_8$ are $\zeta(2)^4, \zeta(5)\zeta(3), \zeta(3)^2\zeta(2), \zeta(6,2)$ and the only 2-3 Lyndon word is
\begin{eqnarray}\label{zeta 6,2}
    \zeta(2,3,3) = \frac{27}{4}\zeta(6,2) - \frac{45}{2}\zeta(5)\zeta(3) + 2\zeta(3)^2\zeta(2) + \frac{1111}{350}\zeta(2)^4.
\end{eqnarray}

In each case the new 2-3 Lyndon word is providing with the new element in each algebra. For $\zeta_9^+$ the generators are $\zeta(7)\zeta(2), \zeta(3)^3, \zeta(5)\zeta(2)^2, \zeta(3)\zeta(2)^3, \zeta(9)$ and the new transcendental element 
$$ \zeta(2,2,2,3) = \frac{641}{16}\zeta(9) - 30\zeta(7)\zeta(2) + \frac{18}{5}\zeta(5)\zeta(2)^2 - \frac{3}{35}\zeta(3)\zeta(2)^3 $$    
    
    And for degree 10 they are $\zeta(8,2), \zeta(3)^2\zeta(2)^2, \zeta(7)\zeta(3), \zeta(5)^2, \zeta(5)\zeta(3)\zeta(2), \zeta(2)^5, \zeta(6,2)\zeta(2)$ and 
    
\begin{eqnarray}\label{zeta 8,2}
\zeta(2,2,3,3) = -\frac{873}{64}\zeta(8,2) + \frac{2037}{32}\zeta(7)\zeta(3) + \frac{1737}{32}\zeta(5)^2 - 24\zeta(5)\zeta(3)\zeta(2) + \frac{3}{5}3/5\zeta(3)^2\zeta(2)^2 - \frac{56643}{7700}\zeta(2)^5
\end{eqnarray}
\subsection{The Broadhurst-Kreimer conjecture}    
	
I used the 2-3 conjecture to find a transcendental basis for the algebra of the MZV as follow.
\begin{center}
    \begin{tabular}{ |c|c|c|} \hline
        Degree & 2-3 Lyndon word & Transcendent basic element \\ \hline
         2; &$\zeta(2)$& $\zeta(2)$\\ \hline
         3; &$\zeta(3)$& $\zeta(3)$\\ \hline
         5; &$\zeta(2,3)$& $\zeta(5)$\\ \hline
         7; &$\zeta(2, 2, 3)$& $\zeta(7)$\\ \hline
         8; &$\zeta(2, 3, 3)$& $\zeta(6,2)$\\ \hline
         9; &$\zeta(2, 2, 2, 3)$& $\zeta(9)$\\ \hline
         10; &$\zeta(2, 2, 3, 3)$& $\zeta(8,2)$\\ \hline
         11; &$\zeta(2, 3, 3, 3), \zeta(2,2,2,2,3)$& $\zeta(8,2,1), \zeta(9,2)$\\ \hline
         12; &$\zeta(2, 2, 2, 3,3), \zeta(2,2,3,2,3)$& $\zeta(8,2,1,1), \zeta(10,2)$\\ \hline
    \end{tabular}
\end{center}
    
\subsection{Conclusion}

I used the \textit{shuffle minus stuffle} product to generate all the known relations of the MZV. Using the xtaylor algorithm I expressed a non-Lyndon word as a product and sum of Lyndon words. Finally I saw that finding the Gr\"obner. basis of degree $n$ for the ideal of the relations is equivalent to solve the linear system given by the relations of degree $n$ without non-Lyndon words (i.e replacing each non-Lyndon words by its factorization) and the Gr\"obner. basis of degree less or equal $n$. All the programs and this paper are available online at:
\begin{center}
http://www.csd.uwo.ca/~gcombar/ResearchMZV/
\end{center}
You can check the Gr\"obner with the relations found by Petitot and Minh in the maple file http://www.csd.uwo.ca/~gcombar/ResearchMZV/mzv16.m. I would like to continue with this research in the future with a more power machine. \\ \\

Germ\'an Combariza.

Computer science department. 

Western University. London ON Canada.

combariza@gmail.com


\end{document}